\def\real{{\tt I\kern-.2em{R}}}
\def\nat{{\tt I\kern-.2em{N}}}
\def\realp#1{{\tt I\kern-.2em{R}}^#1}
\def\natp#1{{\tt I\kern-.2em{N}}^#1}
\def\hyper#1{\ ^*\kern-.2em{#1}}

\def\hyperrealp#1{{\tt ^*{I\kern-.2em{R}}}^#1} 

\def\hypernatp#1{{{^*{{\tt I\kern-.2em{N}}}}}^#1}

\def\leaderfill{\leaders\hbox to 1em{\hss.\hss}\hfill}
\def\srealp#1{{\rm I\kern-.2em{R}}^#1}

\def\power#1{{{\cal P}(#1)}}

\def\qed{{\vrule height6pt width3pt depth2pt}\par\medskip}
\def\pars{\par\smallskip}
\def\parm{\par\medskip}

\def\ref#1{$^{#1}$}

\def\m@th{\mathsurround=0pt}
\def\rightarrowfill{$\m@th \mathord- \mkern-6mu \cleaders\hbox{$\mkern-2mu 
\mathord- \mkern-2mu$}\hfil \mkern-6mu \mathord\rightarrow$}
\def\leftarrowfill{$\mathord\leftarrow
\mkern -6mu \m@th \mathord- \mkern-6mu \cleaders\hbox{$\mkern-2mu 
\mathord- \mkern-2mu$}\hfil $}
\def\noarrowfill{$\m@th \mathord- \mkern-6mu \cleaders\hbox{$\mkern-2mu 
\mathord- \mkern-2mu$}\hfil$}
\def\orgate{$\bigcirc \kern-.80em \lor$}
\def\andgate{$\bigcirc \kern-.80em \land$}
\def\inverter{$\bigcirc \kern-.80em \neg$}

\def\id{\par\hangindent2\parindent\textindent}
\def\textindent#1{\indent\llap{#1}}
\magnification=\magstep1
\tolerance 10000
\baselineskip  12pt
\hoffset=.25in
\hsize 6.00 true in
\vsize 8.85 true in
\font\eightrm=cmr9
\centerline{\bf General Logic-Systems and Consequence Operators}\par\bigskip 
\centerline{Robert A. Herrmann}\parm
\centerline{Mathematics Department}
\centerline{U. S. Naval Academy}
\centerline{572C Holloway Rd.}
\centerline{Annapolis,  MD 21402-5002}
\centerline{1 DEC 2005}\bigskip
{\leftskip=0.5in \rightskip=0.5in \noindent {\eightrm {\it Abstract:} In this paper, general logic-systems are investigated. It is shown that there are infinitely many finite consequence operators defined on a fixed language $L$ that cannot be generated from any finite logic-system. It is shown that any set map $C\colon \power {L} \to \power {L}$ is a finite consequence operator if and only if it is defined by a general logic-system.  
 \par}}\par\bigskip
\noindent{\bf 1. Introduction.}\parm
Usually, two types of consequential operators are investigated, the {\it general} and the {\it finite} ({\it finitary, algebraic}).  In this paper, since the objects will not be embedded formally into a nonstandard structure, roman font will not be used for the informal mathematical symbols. Let $L$ be a nonempty language, $\cal P$ be the power set operator and $\cal F$ the finite power set operator. For a nonempty language $L$, when a general logic-system or a consequence operator is said to be defined on $L$ this means that they are defined on $\power {L}$.\parm
{\bf Definition 1.1.} A mapping $C\colon \power{L} \to \power {L}$ is a general consequence operator (or closure operator) if for each $X,\ Y \in 
\power {L}$\pars
\indent\indent (i) $X \subset C(X) = C(C(X)) \subset L$; and if\pars
\indent\indent (ii) $X \subset Y$, then $C(X) \subset C(Y).$\pars
\noindent A consequence operator C defined on L is said to be {\it finite} ({\it finitary}, or {\it algebraic}) if it satisfies\pars
\indent\indent (iii) $C(X) = \bigcup\{C(A)\mid A \in {\cal F}(X)\}.$\par\medskip
{\bf Remark 1.2.} The above axioms (i), (ii), (iii) are not independent. Indeed, 
(i), (iii) imply (ii). Hence, the finite consequence operators defined on a specific language form a subset of the general operators.\pars
In (Herrmann, 2001), Section 3, logic-systems for a nonempty language are defined by means of a nonempty finite collection of $n$-ary relations $ {\bf RI} = \{R_1,\ldots, R_k\}$, the rules of inference, where each $R_i \subset L^n,\ n \geq 1.$ The set $\bf RI$ includes a possibly nonempty set $R^1,$ (defined as a unary relations, where $(a) = a$), and using this set along with an informally described algorithm deductions are obtain for each $X \in \power {L}.$  It is shown there that each such logic-system generates a finite consequence operator $C$ that yields the same deductions when $C$ is applied to any $X \in \power {L}.$ Since there are but finitely many rules of inference, define such a logic-system as a {\it finite logic-system} and the set of rules $\bf RI$ as the {\it finite rules of inference}. If $\bf RI$ is a finite or infinite set of $n$-ary relations, then $\bf RI$ is called the {\it general rules of inference} and the logic-system that uses such a $\bf RI$ is called a {\it general logic-system.} The symbol $\bf RI$ also denotes the logic-system with its defined processes. As shown in the same manner as for a finite logic-system (Herrmann, 2001), any nonempty general logic-system $\bf RI$ generates (defines) a finite consequence operator the consequences of which are the same as those obtained from $\bf RI.$ 
 \parm
\noindent {\bf 2. General logic-systems.}\parm
 In general, for any $n \geq 1$, $R^n$ denotes a nonempty subset of $L^n.$ I note that the notion of ``effectiveness'' used in Herrmann (2001) can be removed and replaced with either choice or simple conditional statements.

For any consequence operator $C$ generated by a logic-system $\bf RI$, a $R^n \in \bf RI$ {\it applies trivially} to $X\subset L$ if using the algorithm described in Herrmann (2001) restricted to $R^n$ for members of $X$, there is no $a \in L$ such that $a\notin X \cup C(\emptyset).$ For any $X \subset L,$ $R^1$ applies trivially to $X$ since $R^1\subset C(\emptyset).$  Also, if $R^n$ applies trivially to $X\subset L$, then $R^n$ applies trivially to each $Y \subset X.$ \pars
 Consider the natural numbers $\nat$ and a nonempty language $L_N$ each member of which denotes a member of $\nat$. In what follows, the defined standard symbolic forms for members of $\nat$ are used. For $L_N$, let ${\cal S}^f_N$ be the set of all finite rules of inference as defined on $L_N$ and ${\cal S}_N$ the set of all general rules of inference defined on $L_N.$ In all that follows, the superscript denotes the ``arity'' of a relation.  \parm

{\bf Definition 2.1.} The set ${\bf RI}= \{R^n\mid 0<n \in \nat\}$ is defined by induction.\pars
\indent\indent (1) ($n =1$),\  $R^1 = \emptyset$.\pars
\indent\indent (2) ($n>1$),\  $R^n = \{(a_1,\ldots,a_i,\ldots, a_n)\},$ where $a_i = {{n^2-n-4}\over{2}} +i,\ 1\leq i \leq n.$\pars
\parm

{\bf Theorem 2.2.} {\it If $C$ is the finite consequence operator defined on $L_N$ by the infinite general logic-system $\bf RI,$ then there does not exist ${\bf RI_1} \in {\cal S}^f_N$ such that for, $C_1$, its generated consequence operator, $C_1 = C.$}\pars 

Proof. Let ${\bf RI_1} \in {\cal S}^f_N,$ $C_1$ be the finite consequence operator generated by $\bf RI_1$ and $C_1 = C.$  For such a $C_1$, $C_1(\emptyset) = C(\emptyset)=\emptyset$ (i.e. all unary relations are empty.) Since ${\bf RI_1} \in {\cal S}^f_N$ and $C_1(\{0\})= \{0,1\},$ then $\bf RI_1$ contains a nonempty p-ary relation, where $p\geq 2$ is the maximum arity for all of the members of $\bf RI_1.$  Let $k > p$ and $X = \{a_i\mid (a_i = {{k^2 - k -4}\over{2}}+i)\land(i \in \nat)\land(1\leq i\leq k-1)\}.$ Then $C_1(X) = C(X)=\{a_i\mid (a_i = {{k^2 - k -4}\over{2}}+i)\land(i \in \nat)\land(1\leq i\leq k)\}.$  From the definition of $C$, since no member of $X$ appears as any coordinate in any $R^n\in {\bf RI},$ where $n\not= k$ and for $n= k$ they all appear as distinct coordinates, then for each $Y \subset X$ such that $\vert Y\vert = k-2,$ if follows that  for $n \geq 1,$ that $R^n$ applies trivially to $Y$. However, for such $Y$, $C_1(Y) = C(Y)= Y.$  Hence, for every $n \geq 1,$ $R^n \in {\bf RI_1}$ applies trivially to such $Y$. Thus, since $X \subset C_1(X) = C(X)\not= X,$ then there exists one or more $R^q,\ q \geq k,$ that uses all members of $X$ as coordinates and yields a finite step deduction that $C_1(X) = \{a_i\mid (a_i = {{k^2 - k -4}\over{2}}+i)\land(i \in \nat)\land(1\leq i\leq k)\}.$ This contradicts the definition of $p$ and the result follows. \qed   

{\bf Corollary 2.2.1.} {\it There are infinitely many finite consequence operators that cannot be generated by a finite logic-system.}\pars
Proof. Simply change definition 2.1. For each $m \in \nat$, (1) $(n=1),\ R^1_m=\emptyset.$ (2) $(n>1), \  R^n_m =\{(a_1,\ldots,a_i)\}, \ a_i = m + {{n^2 -n-4}\over{2}} +i,\ 1\leq i \leq n.$ \qed

Let $C$ be a general consequence operator defined on $L$. Define $C(\emptyset) =R^1$. Next, for each $\emptyset \not= X \in {\cal F}(L),$ let $\vert X \vert = n.$ Then consider a finite choice (i.e. finite sequence) $X = \{x_1,\ldots, x_n\}.$ Define the corresponding $(n+1)$-ary relation $R^{n+1}_X(C)$ as follows: If $C(X) \not= X\cup R^1$, let $\{(x_1,\ldots,x_n,y)\mid y \in C(X)-(X\cup R^1)\}$ and $\emptyset$ otherwise. Let, for each $n\geq 1,\ n \in {\bf N},$ $R^{n+1}_F(C) = \bigcup\{R^{n+1}_X(C)\mid (\vert X \vert =n)\land(X \in {\cal F}(L))\}.$ Then let ${\bf RI^*} = \{R^1\} \cup \{R^{n+1}_F(C) \mid (n \geq 1)\land (n\in {\bf N})\}.$ Obviously, $ \bf RI^*$ is not unique.  (Note: This definition is somewhat similar to the definition given by {\L}os and Suszko for general consequence operators on an organized sentential language. But, it does not call for any additional closure conditions.)\parm

{\bf Theorem 2.3.} {\it Let $C$ be a finite consequence operator defined on nonempty $L$  and $\bf RI^*$ the general rules of reference as defined above by $C.$ If $C^*$ is the finite consequence operator generated by $\bf RI^*$, then $C^* = C.$}\pars 
Proof.  Recall that for two consequence operators, $C',C''$  defined on $L,$ $C'\leq C''$ if and only if for each $X \in {\cal P}(L),\ C'(X) \subset C''(X).$ Indeed, $\langle{\cal C}_f(L),\leq\rangle$ is a sublattice of the lattice $\langle{\cal C}(L),\leq\rangle$. For a language $L$, let $\bf RI^*_1$ be the rules of inference defined by a consequence operator $C_1$, where $C_1$ is not necessarily finite. The following is established by induction on the number of steps in a deduction.\parm

Consider $\bf RI^*_1$, where $C_1$ is not necessarily finite. Let $X \in {\cal P}(L)$ and $\{b_1,\ldots,b_n\}$, where all members are always assumed distinct, be an $\bf RI^*_1$-deduction from $X$.  Then $b_i \in C_1(X)$ for each $i$ such that $1\leq i \leq n.$ \parm

(1) Let $n =1.$ Then $b_1 \in X$ or $b_1 \in R^1.$ Since, by insertion,  $X \cup R^1 \subset C_1(X),$ then $b_1 \in C_1(X).$ \pars

(2) Consider a deduction $\{b_1,\cdots,b_{n+1}\}$ from $X$ and assume the strong induction hypothesis that $b_i \in C_1(X)$ for each $i$ such that $1 \leq i \leq n.$ Then either $b_{n+1} \in X \cup R^1$ or not. If $b_{n+1} \in X \cup R^1$, then, as in (1), $b_{n+1} \in C_1(X).$ Otherwise, by definition of logic-system deduction, there exists some $R^{k+1}_F(C_1) \in  {\bf RI^*_1}$ and a $(y_1,\ldots,y_{k+1}) \in R^{k+1}_F(C_1)$ such that nonempty $M = \{y_1,\ldots,y_k\} \subset \{b_1\ldots,b_n\} \cup X \cup R^1$ and $b_{n+1} \in C_1(M).$ However, $\{b_1,\ldots,b_n\} \cup X \cup R^1 \subset C_1(X)$ yields that $C_1(M)\subset C_1(\{b_1\ldots,b_n\} \cup X \cup R^1) \subset C_1(C_1(X))= C_1(X).$ Hence, $b_{n+1} \in C_1(X)$ and the result follows by induction.\pars

Now let $C^*$ be the finite consequence operator generated by the $C$ generated $\bf RI^*,$ $X \in {\cal P}(L)$, and $x \in C^*(X)$. Then, from the definition of $C^*$, there is an $ \bf RI^*$-deduction $\{b_1,\ldots,b_n\}$ from $X$ such that $b_n = x.$ Hence, $x \in C(X).$
Thus, $C^* \leq C.$ \pars

Conversely, for finite $C$, let $x \in C(X)$. If $x \in X \cup R^1$, where $R^1 = C(\emptyset)$, then $x \in C^*(X)$ by insertion. Hence, assume that $x \notin X \cup R^1.$ Since $C$ is finite, there is some finite $F \subset X,$ of smallest cardinality, such that $x \in C(F)$. The set $F\not=\emptyset,$ since $C(\emptyset) = R^1.$ Therefore, $m = \vert F \vert \geq 1.$ Thus, there is an $r \in R^{m+1}_F$ such that $p_i(r)\in F,\ 1\leq i\leq m$ and $p_{m+1}(r) = x$ and from the definition of $C^*,$ $x \in C^*(F) \subset C^*(X).$ Thus, $C \leq C^*$. Hence, $C^* = C$ and the proof is complete.\parm 

\noindent {\bf 3. Generating logic-systems.}\parm

In the physical sciences, the set $\bf RI$ is usually not defined explicitly. In actual practice, a physical argument simply claims that a specific finite set of statements - the conclusions - (among other names) is ``deduced'' from another finite set of statements $S$. These statements $S$ can contain members from a fixed set of statements $A$ - a unary relation - where $A$ can be further partitioned. The actual hypotheses $H= S-A.$ (In formal logic, the $A$ is the entire set of statements generated by the axioms.) The actual rules of the logic are not usually stated. It is assumed that after refinements via peer evaluation that the vast majority of the members of a specific science-community would accept the ``derivation.'' It is obvious how one would construct a general $\bf RI$ from collections of such derivations. For example, let $\{h_1,\ldots,h_n\}\cup \{a_{n+1},\ldots,a_{n+k}\} = S,\ n\geq 1,\ k\geq 1,$ where $H = \{h_1,\ldots, h_n\},\ \{a_{n+1},\ldots,a_{n+k}\} \subset A$ and $\{b_1,\ldots,b_m\}=B$ is non-trivially deduced from $S$. Then one can construct a $R^{n+k+1}\in \bf RI$ that contains $\{(h_1,\ldots,h_n,a_{n+1},\ldots,a_{n+k},x)\mid x \in B\}.$ However, under our definition of how logic-systems are employed for deduction, many other relations can also generate each member of $B$. For example, consider $\{(h_1,\ldots,h_n,a_{n+1},\ldots,a_{n+k},b_1)\}$ and $\{(b_1,x)\mid x \in B\}.$Since there must be at the least one member in $S$ for there to be any non-trivial deduction, then, for certain $S,$  it might be discovered that the set $B$ can be deduced from the set $\{h_1\}\in H$ and $\{a_1\}\in A$. Hence, from the definition of logic-systems and how they are used to generate deductions, you could also have a 3-ary relation $R^3 \in \bf RI$ such that $\{(h_1,a_1,x)\mid x \in B\} \subset R^3.$ Consequently, in general for non-trivial deduction from premises, it is merely assumed that for any set $X \in {\cal P}(L)$, if there is nonempty $B\in {\cal P}(L),\ B\cap X =\emptyset$ and each member of $B$ is claimed to be deduced from finitely many members of $X,$ using finitely members of a auxiliary set $A,$ then, at the least, there are $n$-ary relations in $\bf RI$ that contain members that generate the members of $B$ according to the algorithm stated in Herrmann (2001). \pars
Theorem 2.3 yields an obvious question relative to theorem 3.8 in Herrmann (2004). Let $\{{\bf RI_i}\}$ be a nonempty collection of rules of inference, where each $\bf RI_i$ is defined on a nonempty language $L_i.$ Let each $\bf RI_i$ generate a corresponding finite consequence operator $C_i$. If $C$ is the finite consequence operator on $\bigcup \{L_i\}$ generated by the rules of inference $\bigcup \{\bf RI_i\}$, is $C = \bigvee_W \{C_i\}$?    
\parm

\centerline{\bf References}\parm
\id{H}errmann, Robert A. (2004), ``The best possible unification for any collection of physical theories,'' {\it Internat. J. Math. and Math. Sci.}, 17(2004):861-721.\hfil\break
http://www.arxiv.org/abs/physics/0306147 http://www.arxiv.org/abs/physics/0205073\smallskip 

\id{H}errmann, Robert A. (2001), ``Hyperfinite and Standard Unifications for Physical Theories,'' {\it International Journal of Mathematics and Mathematical Sciences}, 28(2):93-102. 
http://www.arXiv.org/abs/physics/0105012\smallskip

\id{\L}os, I. and R. Suszko, (1958), ``Remarks on sentential logic,'' Indagationes Mathematical, 20:177-183
\end